\documentclass[10pt]{article}
\usepackage{mathrsfs}
\usepackage{amssymb,amsfonts,amsmath, psfrag,eepic,colordvi,epsfig}
\topskip 
\numberwithin{equation}{section}
\newtheorem{theorem}{Theorem}[section]

\newtheorem{conjecture}[theorem]{Conjecture}

\newtheorem{lemma}[theorem]{Lemma}

\begin{document}
\parskip 7pt

\pagenumbering{arabic}
\def\sof{\hfill\rule{2mm}{2mm}}
\def\ls{\leq}
\def\gs{\geq}
\def\SS{\mathcal S}
\def\qq{{\bold q}}
\def\MM{\mathcal M}
\def\TT{\mathcal T}
\def\EE{\mathcal E}
\def\lsp{\mbox{lsp}}
\def\rsp{\mbox{rsp}}
\def\pf{\noindent {\it Proof.} }
\def\mp{\mbox{pyramid}}
\def\mb{\mbox{block}}
\def\mc{\mbox{cross}}
\def\qed{\hfill \rule{4pt}{7pt}}
\def\pf{\noindent {\it Proof.} }
\textheight=22cm

\begin{center}
{\Large\bf On the log-concavity
 of $n$-th root of a sequence  }
\end{center}

\begin{center}

Ernest X. W. Xia

Zuo-Ru Zhang\\[8pt]
School of Mathematical Science\\
Hebei Normal University\\
 Shijiazhuang 050024, P. R. China\\[6pt]

{\tt zrzhang@hebtu.edu.cn}

\end{center}

\noindent {\bf Abstract.}
  In recent years,
   the log-concavity
    of  $\{\sqrt[n]{S_n}\}_{n\geq 1}$
      have been received a lot of
       attention. Very recently,
    Sun posed the following conjecture
     in his new book:
   the sequences $\{
   \sqrt[n]{a_n}\}_{n\geq 2}$
    and $\{  \sqrt[n]{b_n}\}_{n\geq 1}$
     are log-concave,
      where
     \[
a_n:= \frac{1}{n}\sum_{k=0}^{n-1}
 \frac{{n-1\choose k}^2{n+k\choose k}^2
  }{4k^2-1}
     \]
and
     \[
b_n:= \frac{1}{n^3}\sum_{k=0}^{n-1}
  (3k^2+3k+1){n-1\choose k}^2
  {n+k\choose k}^2.
     \]
In this paper, two methods, semi-automatic and
  analytic  methods, are used to confirm   Sun's conjecture. The semi-automatic
  method relies on a      criterion on the log-concavity   of  $\{\sqrt[n]{S_n}\}_{n\geq 1}$
   given by us and  a mathematica package  due to Hou and Zhang, while the
analytic  method relies on a result due to Xia.

\noindent {\bf Keywords}: log-concavity,
 inequality, combinatorial sequence.

\noindent { 2010 Mathematics subject
 classification}:   05A20;
05A10.

\section{Introduction}

\allowdisplaybreaks

Recall that a positive sequence
 $\{S_n\}_{n\geq 1}$
   is log-concave if for $n\geq 1$,
   \begin{align}\label{1-1}
\frac{S_{n+1}}{S_{n}}\geq
 \frac{S_{n+2}}{S_{n+1}}.
   \end{align}
Meanwhile, the sequence
 $\{S_n\}_{n\geq 1}$
 is called strictly
log-concave   if the inequality in
  \eqref{1-1}
   is strict for all  $n\geq 1$.

In 1982, Firoozbakht conjectured
 that for $n\geq 1$,
  \[
\sqrt[n]{p_n}>\sqrt[n+1]{p_{n+1}},
  \]
  where $p_n$ denotes  the $n$th
   prime number.  In 2013,
     Sun \cite{Sun-0}
     proved that the sequence
      $\{\sqrt[n]{\alpha_n}\}_{n\geq 2}$
        is strictly decreasing
         and the sequence
      $\{\sqrt[n+1]{\alpha_{n+1}}/
      \sqrt[n]{\alpha_n}\}_{n\geq 5}$
       is strictly increasing, where
       $\alpha_n=p_1+p_2+\cdots +p_n$.
        Motivated by those works,
         Sun \cite{Sun-1} posed a number
          of conjectures on
           monotonicity of the sequences
            $\{ \sqrt[n]{S_n}\}_{n\geq 1}$
            and $\{ \sqrt[n+1]{S_{n+1}
             }/\sqrt[n]{S_n}
              \}_{n\geq 1}$,
               where $\{S_n\}_{n\geq 1}$
                are  number-theoretic
                or
combinatorial sequences
 of positive integers. Some of those conjectures
  were confirmed by
    Chen,
   Guo and Wang \cite{Chen},
    Hou, Sun and Wen \cite{Hou},
    Luca  and St\u{a}nic\u{a} \cite{Luca},
      Wang and Zhu
    \cite{Wang}, and Xia \cite{Xia}.

 Recently, Sun posed the following conjecture on two sequences
    $\{a_n\}_{n\geq 1}$  and  $\{b_n\}_{n\geq 1}$ \cite{Sloane}
in his new book \cite{Sun-2}.

\begin{conjecture} \label{conjecture}
 Define
\[
a_n:= \frac{1}{n}\sum_{k=0}^{n-1}
 \frac{{n-1\choose k}^2{n+k\choose k}^2
  }{4k^2-1}
     \]
and
     \[
b_n:= \frac{1}{n^3}\sum_{k=0}^{n-1}
  (3k^2+3k+1){n-1\choose k}^2
  {n+k\choose k}^2.
     \]
Then the sequences $\{a_{n+1}/a_n\}_{n\geq 3
 }$
  and  $\{b_{n+1}/b_n\}_{n\geq 1}$
  are both strictly increasing to the limit
   $17+12\sqrt{2}$, while the sequences
 $\{\sqrt[n+1]{a_{n+1}}
 /\sqrt[n]{a_n}\}_{n\geq 2
  }$
  and  $\{\sqrt[n+1]{b_{n+1}}/
   \sqrt[n]{b_n}\}_{n\geq 1}$
   are strictly    decreasing to 1.

\end{conjecture}

The aim of this paper is to prove
 Sun's conjecture
  by using two methods:  semi-automatic
   and
  analytic
  methods. The semi-automatic
   proof require
    a
      criterion on the log-concavity
       of  $\{\sqrt[n]{S_n}\}_{n\geq 1}$
 and
  a mathematica package   $\mathbf{p.rec}$
   due to Hou and Zhang \cite{Hou-Zhang},
    while the  analytic proof
     require a result
      proved by Xia \cite{Xia}.

The rest of this paper is organized as follows.
In Section 2, we
 prove a
      criterion on the log-concavity
       of  $\{\sqrt[n]{S_n}\}_{n\geq 1}$.
 In Section 3, we
 present a semi-automatic
   proof  of   Conjecture \ref{conjecture} by combining
  the
    criterion given in Section 2
     and  a mathematica package   {\bf
      p.rec}
   given by  Hou and Zhang \cite{Hou-Zhang}.
     In Section 4,  we present
     an   analytic proof
      of Conjecture \ref{conjecture}
       by using a result given by Xia
        \cite{Xia}.

 \section{A
      criterion on the log-concavity
       of $n$-th root of a sequence}

In this section, we pose  the following
 criterion on the log-concavity
       of $\{\sqrt[n]{S_n}\}_{n\geq 1}$
        which will be used to
         confirm Sun's conjecture.

 \begin{theorem}\label{Th-1}
  Let $\{S_n\}_{n\geq 0}$
   be a sequence of positive integers.
  Suppose that  there exist
    positive numbers $a_0,  k_0, N_0,m$
     and  a function
   \[
v(n,m):=v(n)=\frac{c_1}{n}+\cdots +\frac{c_m}{ n^m}
   \]
    such that
      for $n\geq N_0$,

 (1) $a_0+v(n)-\frac{1}{n^m}<
 \frac{S_{n+1}}{
 S_n}<a_0+v(n)+\frac{1}{n^m}
 <a_0\frac{n^{k_0}}{
  (n+1)^{k_0}}$;

 (2) $-\frac{1}{2}<v(n)-\frac{1}{n^m}$;

 (3) $S_{N_0}\leq a_0^{N_0}$;

 (4) $\frac{v(n)+\frac{1}{n^m}}{
  a_0(n+1)}-
  \frac{v(n-1)-\frac{1}{(n-1)^m}}{
  a_0(n-1)}+ \frac{(v(n-1)-\frac{1}{(n-1)^m})^2
   }{
  a_0^2(n-1)}-\frac{2k_0\ln(n)}{n(n^2-1)}<0$.
\\
Then   $\{\sqrt[n]{S_n}\}_{n\geq N_0}$
 is log-concave and
 \[
\lim_{n\rightarrow \infty}
 \frac{\sqrt[n+1]{S_{n+1}}}{\sqrt[n]{S_n}}=1.
\]

 \end{theorem}

\noindent{\it Proof.}
  Since $\ln(1+x)\leq x$  for $x>-1$, we have
\begin{align}
\ln\left(a_0+v(n)+\frac{1}{n^m}\right)
 =&
\ln(a_0)+\ln\left(1+
\frac{v(n) +\frac{1}{n^m}}{a_0}\right) \nonumber
\\[6pt]
\leq & \ln(a_0)+\frac{v(n)
 +\frac{1}{n^m}}{a_0}. \label{2-1}
\end{align}
On the other hand, it is easy to prove
 that
  for  $x> -\frac{1}{2}$,
   \[\ln (1+x)\geq x-x^2 .\]
Thus,
\begin{align}
\ln\left(a_0+v(n)-\frac{1}{n^m}\right)
 =& \ln(a_0)+\ln\left(1+
\frac{v(n) -\frac{1}{n^m}}{a_0}\right)
\nonumber
\\[6pt]
> & \ln(a_0)+\frac{v(n) -
\frac{1}{n^m}}{a_0}-\left(\frac{v(n) - \frac{1}{n^m}}{a_0}
 \right)^2.  \quad ({\rm by\ Condition
  \ }\ (2)) \label{2-2}
\end{align}
In view of Conditions (1) and (3),
\begin{align}
S_n=S_{N_0} \prod_{i=N_0}^{n-1}
 \frac{S_{i+1}
  }{S_i}<S_{N_0} \prod_{i=N_0}^{n-1} \left(
  a_0
  \frac{i^{k_0}}{(i+1)^{k_0}}\right)
  =\frac{S_{N_0} a_0^n}{a_0^{N_0} n^{k_0}}
  <\frac{  a_0^n}{  n^{k_0}}.
 \label{2-3}
\end{align}
It is easy to verify that
 for $n\geq N_0$,
\begin{align*}
&\ln\left(
 \frac{\sqrt[n+1]{S_{n+1}}\sqrt[n-
 1]{S_{n-1}}}{(\sqrt[n]{S_{n}})^2}\right)
 \nonumber\\[6pt]
 =
&\frac{\ln S_{n+1}}{n+1}
 + \frac{\ln S_{n-1}}{n-1}
 -2\frac{\ln S_n}{n}
 \nonumber\\[6pt]
 =&\frac{\ln(S_{n+1}/S_n)}{
 n+1}-\frac{\ln(S_{n}/S_{n-1})}{n-1}
 +2\frac{\ln(S_n)}{n(n^2-1)}
 \nonumber\\[6pt]
 <&\frac{1}{n+1}
 \ln\left(a_0+v(n)+\frac{1}{n^m}\right)
 -\frac{1}{n-1}
  \ln\left(a_0+v(n-1)
  -\frac{1}{(n-1)^m}\right)
  \nonumber\\[6pt]
  &\qquad
  +2\frac{\ln(S_n)}{n(n^2-1)}
  \quad ({\rm by\ Condition \ }\ (1))
 \nonumber\\[6pt]
 <& \frac{1}{n+1}\left(\ln(a_0)+\frac{v(n)
 +\frac{1}{n^m}}{a_0}\right)
 \nonumber\\[6pt]
 &-\frac{1}{n-1}
 \left(\ln(a_0)+\frac{v(n-1) -
\frac{1}{(n-1)^m}}{a_0} -\left(\frac{v(n-1 ) -
\frac{1}{(n-1)^m}}{a_0}
 \right)^2\right)
 \nonumber\\[6pt]
& +2\frac{\ln(a_0^n/n^{k_0})}{n(n^2-1)}
 \quad ({\rm by\  }\ \eqref{2-1}-\eqref{2-3})\nonumber\\[6pt]
  =&\frac{v(n)+\frac{1}{n^m}}{
  a_0(n+1)}-
  \frac{v(n-1)-\frac{1}{(n-1)^m}}{
  a_0(n-1)}+ \frac{(v(n-1)-\frac{1}{(n-1)^m})^2
   }{
  a_0^2(n-1)}-\frac{2k_0\ln(n)}{n(n^2-1)}
  \nonumber\\[6pt]
  <&0, \quad ({\rm by\ Condition \ }\ (4))
\end{align*}
which implies that
 for $n\geq N_0$,
   \[
\frac{\sqrt[n+1]{S_{n+1}}
 }{\sqrt[n]{S_{n}}}<
 \frac{\sqrt[n]{S_{n}}
  }{\sqrt[n-1]{S_{n-1}}}.
   \]
Thus $\{S_n\}_{n\geq N_0}$
 is log-concave.

%

 Note that for $n\geq N_0$,
 \begin{align}
\ln\frac{\sqrt[n+1]{S_{n+1}}
 }{ \sqrt[n]{S_n}} &=\frac{1}{n+1} \ln
(S_{n+1}/S_n)-\frac{1}{n(n+1)}\ln (S_n) \nonumber\\[6pt]
 &>\frac{1}{n+1} \ln
(S_{n+1}/S_n)-\frac{1}{n(n+1)}\ln
 \left(\frac{  a_0^n}{  n^{k_0}}
 \right) .\qquad ({\rm by}
 \ \eqref{2-3}) \label{2-4}
 \end{align}
In addition,
\begin{align}
\ln\frac{\sqrt[n+1]{S_{n+1}}
 }{ \sqrt[n]{S_n}} <\frac{1}{n+1} \ln
(S_{n+1}/S_n)  \label{2-5}
\end{align}
since $S_n\geq 1$. Note that
\begin{align}\label{2-6}
\lim_{n\rightarrow \infty} \frac{1}{n+1} \ln (S_{n+1}/S_n)=
\lim_{n\rightarrow
 \infty }\frac{1}{n+1} \ln
(S _{n+1}/S_n)-\frac{1}{n(n+1)}\ln
 \left(\frac{ a_0^n}{n^{k_0}}
  \right)=0.
\end{align}
It follows from \eqref{2-4}-\eqref{2-6} that
\[
\lim_{n\rightarrow \infty}\ln\frac{\sqrt[n+1]{a_{n+1}}
 }{ \sqrt[n]{a_n}}=0
\]
and
\[
\lim_{n\rightarrow
 \infty} \frac{\sqrt[n+1]{a_{n+1}}
 }{ \sqrt[n]{a_n}}=1.
\]
 This completes the
  proof. \qed

\section{A semi-automatic
   proof  of   Sun's conjecture}

In this section, we present a semi-automatic proof  of   Sun's conjecture based on Theorem \ref{Th-1}
 as well as a mathematica
 package {\bf p.rec} given by Hou and Zhang \cite{Hou-Zhang}.

\subsection{The  mathematica
 package p.rec}

%

Recall that  a sequence
  $\{T_n\}_{n\geq 0}$
    is called asymptotically $r
     $-log-convex if
     $\mathscr{L}^k(\{T_n\}_{n\geq N})$
are non-negative sequences for
 $1 \leq
   k \leq  r$ and a certain integer $N$,
    where $\mathscr{L}$ is an operator
     defined by
     \[
\mathscr{L}(\{T_n\}_{n\geq 0}=\{T_{n+2}T_n
 -T_{n+1}^2\}_{n\geq 0}
 \qquad {\rm and}\ \qquad
   \mathscr{L}^k(\{T_n\}_{n\geq 0})
   = \mathscr{L}(
     \mathscr{L}^{k-1}(\{T_n\}_{n\geq 0})).
        \]
In \cite{Hou-Zhang} Hou and Zhang
 give a general method to prove
 the asymptotic $r$-log-convexity
 of a sequence.
 As an application, they
   considered
  the problem of proving the
  asymptotic $r$-log-convexity
  of P-recursive sequences
   $\{S_n\}_{n\geq0}$ defined by
\begin{equation}\label{prec}
S_{n+t}=r_0(n)S_n+r_1(n)S_{n+1}
 +\cdots+r_{t-1}(n)S_{n+t-1},
\end{equation}
where   $r_i(n)$  are rational functions of $n$.
More precisely,   Hou and Zhang \cite{Hou-Zhang}
 presented a method for   proving the asymptotic   $r$-log-convexity
 of P-recursive sequences  based on   the asymptotic expansions of the
P-recursive sequences given by Birkhoff and Trjitzinsky
\cite{Birkhoff} and developed by Wimp and Zeilberger \cite{Wimp}.
 As an example, they confirmed a part of a conjecture
   posed by Chen and Xia \cite{Chen}.
   They also posed     an algorithm to compute an
       explicit $N_0 = N_0(r)$
   such that $\{S_n\}_{n\geq N_0}$
    is $r$-log-convex
     for   bound preserving sequence
      $\{S_n\}_{n\geq 0}$.
       Define
       \[
 f_n: = \sum_{i=m_0}^m \frac{c_i}{n^{i/\rho}} - \frac{1}{n^{m/\rho}}, \quad
 g_n: = \sum_{i=m_0}^m \frac{c_i}{n^{i/\rho}} + \frac{1}{n^{m/\rho}},
\]
where $\rho$ is   a positive integer,
  $m_0$ is an integer and $c_i$ are real
numbers.
 Recall that the sequence
  $\{S_n\}_{n \ge 0}$ defined by \eqref{prec}
    is called
    a \emph{bound preserving} sequence if for $m$ sufficiently large,
     there exists an integer $N_0$ such that
\[
\frac{r_0(n)}{u^{(0)}_n u^{(0)}_{n+1} \cdots u^{(0)}_{n+t-2}}+
\frac{r_1(n)}{u^{(1)}_{n+1} \cdots u^{(1)}_{n+t-2}} +
\cdots+r_{t-1}(n) \geq f_{n+t-1}^{(0)}, \quad \forall\, n \ge N_0,
\]
and
\[
\frac{r_0(n)}{v^{(0)}_n v^{(0)}_{n+1}
\cdots v^{(0)}_{n+t-2}}+
\frac{r_1(n)}{v^{(1)}_{n+1} \cdots v^{(1)}_{n+t-2}} +
\cdots+r_{t-1}(n) \leq g_{n+t-1}^{(0)}, \quad \forall\, n \ge N_0,
\]
where
\[
u^{(i)}_n=
\begin{cases}
g_n, & \mbox{if $r_i(n)>0$ for
$n\geq N_0$,}\\[5pt]
f_n, &  \mbox{if $r_i(n)<0$ for $n \geq N_0 $, }
\end{cases}
\]
and
\[
v^{(i)}_n=\begin{cases}
f_n, & \mbox{if $r_i(n)>0$ for $n\geq N_0$,
  }\\[5pt]
g_n, & \mbox{if $r_i(n)<0$ for $n\geq N_0$.  }
\end{cases}
\]
Hou and Zhang \cite{Hou}
 also
 implemented a Mathematica
  package {\bf p.rec}
  which is accessible at
  \begin{center}
  http://cam.tju. edu.cn/~hou/preprints.html
\end{center}
 to find the above integer $N_0$.    In fact, we can use some
 commands of the package to find    the upper and lower bounds of
    $S_{n+1}/S_n$.
Such as the command
\begin{verbatim}
        rLogBound[L, n, N, ini_val, r, t]
\end{verbatim}
where
\[
L=N^d-(R_0(n)+R_1(n)N+\cdots+R_{d-1}N^{d-1})
\]
is the annihilating operator for the P-recursive sequence corresponds to the recurrence relation. Here \verb|N| is the forward shift operator of $n$,   \verb|ini_val| is a list of initial values of the sequences and \verb|t| is the number of terms of the asymptotic ratio of the sequence.

\subsection{Upper and lower bounds
 for $a_{n+1}/a_n$
  and $b_{n+1}/b_n$ }

In this subsection, we use two
   mathematica packages {\bf p.rec}
    and   {\bf fastzeil} to deduce
     upper and lower
 bounds for $a_{n+1}/a_n$
  and $b_{n+1}/b_n$.

Since both the two sequences are summations of two proper hypergeometric terms of $n$ and $k$, according to Sister Celine's existence theorem (Theorem 4.4.1 \cite{w-z}), there exists a $k-$ free recurrence of both sequences.  Through Zeilberger's algorithm and a improved  mathematica package
 {\bf fastzeil} given by Paule and Schorn\cite{paule} which is accessible at
 \begin{center}
   https://risc.jku.at/sw/fastzeil/
 \end{center}
 We get the recurrence relations
 for $\{a_n\}_{n\geq 0}
 $ and $\{b_n\}_{n\geq 0}$.
 Using the following command
\begin{align*}
  {\rm In}[1]:=&{\rm zb}[({\rm
   Binomial}[n-1,k]{\rm Binomial}[n+k,k])^2/
   (n(4k^2-1)),{k,0,n-1},n,3];
     \\[6pt]
     {\rm out}[1]=&\{n^3(1+n)(5+2n){\rm
      SUM}[n]
    -(1+n)(5+2n)(62+191n+152n^2+35n^3)
        {\rm SUM}[n+1]
   \nonumber\\[6pt]
     &+(2+n)(1+2n)(88+224n+163n^2+35n^3)
     {\rm SUM}[n+2]
     \nonumber\\[6pt]
     &
   -(2+n)(3+n)^3(1+2n){\rm SUM}[n+3]=0\}
\end{align*}
    The output implies that
    $a_n$ satisfies
  the following recurrence
    relation:
\begin{align}\label{3-3}
a_n=v_1(n)a_{n-1}+v_2(n)a_{n-2}+v_3(n)a_{n-3},
\end{align}
where
\begin{align*}
v_1(n)=&\frac{35n^3-152n^2+191n-62}{ n^3},
\\[6pt]
v_2(n)=&-\frac{ (n-2)(2n-1)(-88+224n
-163n^2+35n^3)}{(n-1)n^3(2n-5)},
\\[6pt]
v_3(n)=&\frac{ (n-3)^3(n-2)(2n-1)}{
 (n-1)n^3(2n-5)}.
\end{align*}
Now, we can use the command {\bf
 rLogBound} to deduce the upper and lower
  bounds for $a_{n+1}/a_n$.
   This is done automatically by
\begin{align*}
{\rm In}[2]:=&{\rm L1}=n^3(2n+5)(n+1)-(2n+5)(n+1) (35n^3 +152n^2
+191n+62)N\nonumber\\[6pt]
 &  +(n+2)(2n+1)(35n^3+163n^2+224n+88)N^2
   -(n+2)(2n+1)(n+3)^3N^3;\\[6pt]
   {\rm In}[3]:=&{\rm rLogBound}[L1, n, N, {0, -1, 1, 9, 61}, 1, 3]\\[6pt]
{\rm Out}[3]=&
 17+12\sqrt{2}-\frac{\frac{153}{2}+54\sqrt{2}
 }{n}+\frac{\frac{675}{4}+\frac{7617\sqrt{2}}{64}
  }{n^2}-\frac{1}{n^2}\nonumber\\[6pt]
  &\leq \frac{a_{n+1}}{a_n}
 \leq 17+12\sqrt{2}-\frac{\frac{153}{2}+54\sqrt{2}
 }{n}+\frac{\frac{675}{4}+\frac{7617\sqrt{2}}{64}
  }{n^2}-\frac{1}{n^2} \quad for \ n\geq 5
  \\[6pt]
  &a_n \ {\rm   preserves\  the \
   bounds \ for} \ n\geq 605\\
  &{\rm the \ bounds \ hold \ for} \
    n\geq 607
  \\
  &\{{\rm True}, \ 607\}
\end{align*}
Therefore, we arrive at
 the following lemma.

\begin{lemma} \label{L-1}
For $n\geq 607$, we have
\begin{align}\label{3-1}
h_1(n)-\frac{1}{n^2} < \frac{a_{n+1} }{a_n
 }<h_1(n)+\frac{1}{n^2},
\end{align}
where
\begin{align}\label{3-2}
h_1(n)=17+12\sqrt{2}-\frac{\frac{153}{2}+54\sqrt{2}
 }{n}+\frac{\frac{675}{4}+\frac{7617\sqrt{2}}{64}
  }{n^2}.
\end{align}
\end{lemma}

Similarly, we can also find the recurrence
 relation for $b_n$.
\begin{align*}
  {\rm In}[4]:=&{\rm zb}[(3k^2+3k+1)({\rm
   Binomial}[n-1,k]{\rm Binomial}
    [n+k,k])^2/n^3,{k,0,n-1},n,3];
    \\[6pt]
    {\rm
     out}[4]=&\{
    (1+n)(5+2n)(11+12n+3n^2)(25+24n+6n^2){\rm SUM}[n]
    \nonumber\\[6pt]
   &-(1+n)(5+2n)(3076+21646n+59512n^2+82777n^3+64134n^4
   +28137n^5 \\[6pt] &+6552n^6+630n^7){\rm SUM}[n+1]
     +(2+n)(1+2n)(5072+30640n+73445n^2
   \nonumber\\[6pt]
   & +93469n^3+68751n^4+29271n^5
    +6678n^6+630n^7){\rm SUM}[n+2]
      \nonumber\\[6pt]
    &
    -(2+n)(3+n)^3(1+2n)(2+6n+3n^2)
     (7+12n+6n^2){\rm
      SUM}[n+3]=0\}
    \end{align*}
The output means that   $b_n$ satisfies
 the following recurrence relation:
\begin{align}\label{3-11}
b_n=u_1(n)b_{n-1}+u_2(n)b_{n-2}
 +u_3(n)b_{n-3},
\end{align}
where
\begin{align*}
u_1(n)=& \frac{(630n^7-6552n^6
 +28137n^5-64134n^4
  +82777n^3-59512n^2+21646n-3076)}{
   (3n^2-12n+11)(6n^2-24n+25)n^3},
   \\[6pt]
   u_2(n)=&   -\frac{
    (n-2)(2n-1) w(n)}{
     (n-1)(2n-5)(3n^2-12n+11)
      (6n^2-24n+25)n^3},
      \\[6pt]
      u_3(n)=& \frac{(n-2)(2n-1)
       (3n^2-6n+2)(6n^2-12n+7)(n-3)^3}{
        (n-1)(2n-5)(3n^2-12n+11)
         (6n^2-24n+25)n^3}
\end{align*}
with
\[
w(n)=630n^7-6678n^6
    +29271n^5-68751n^4+93469n^3
    -73445n^2+30640n-5072.
\]
With the next command we  find  the upper and lower
  bounds for $b_{n+1}/b_n$:
  \begin{align*}
{\rm In}[5]:=&L2=(2n+5)(n+1)(6n^2+24n+25) (3n^2+12n+11)n^3  -
 (630n^7+6552n^6\nonumber\\[6pt]
&+28137n^5+64134n^4 +82777n^3+59512n^2+21646n+3076)
 (2n+5)(n+1)N\nonumber\\[6pt]
 &  +(n+2)(2n+1)(630n^7+6678n^6+29271n^5+68751n^4
 +93469n^3+73445n^2\nonumber\\[6pt]
&+30640n+5072)N^2
    -(n+2)(2n+1)(n+3)^3(6n^2+12n+7)(3n^2+6n+2)N^3;\\[6pt]
   In[6]:=&{\rm rLogBound}[L2, n, N, {1,8,87,1334,25045}, 1, 4]\\[6pt]
{\rm Out}[6]=&
 17+12\sqrt{2}
 -\frac{\frac{85}{2}+
  30\sqrt{2}
  }{n}+\frac{\frac{275}{4}+\frac{3105}{
  32\sqrt{2}}
 }{n^2} -\frac{1}{n^2}\nonumber\\[6pt]
  &\leq \frac{b_{n+1}}{b_n}
 \leq 17+12\sqrt{2}
 -\frac{\frac{85}{2}+
  30\sqrt{2}
  }{n}+\frac{\frac{275}{4}+\frac{3105}{
  32\sqrt{2}}
 }{n^2} +\frac{1}{n^2} \quad for \ n\geq 3
  \\[6pt]
  &b_n \  {\rm preserves\  the \
   bounds \ for }\ n\geq 208\\
  &{\rm the \ bounds \ hold \ for\ }
    n\geq 210
  \\
  &\{{\rm True}, \ 210\}
\end{align*}

Thus, we obtain the following lemma:

\begin{lemma}\label{L-2}
 For $n\geq 210$,
 \begin{align}\label{3-12}
h_2(n)-\frac{1}{n^2} < \frac{ b_{n+1} }{b_n
 }<h_2(n)+\frac{1}{n^2},
\end{align}
where
\begin{align}\label{3-13}
h_2(n)=17+12\sqrt{2} -\frac{\frac{85}{2}+30
 \sqrt{2}
 }{n}+\frac{\frac{275}{4}+\frac{3105
  \sqrt{2}}{64}
  }{n^2}.
\end{align}
\end{lemma}

\subsection{A proof of Conjecture
 \ref{conjecture} }

Now, we are ready to prove Conjecture
 \ref{conjecture} based on Theorem \ref{Th-1}
  and Lemmas \ref{L-1}
   and \ref{L-2}.

By \eqref{3-1} and \eqref{3-12}, we see that
\begin{align*}
\lim_{n\rightarrow \infty}
 \frac{a_{n+1}}{a_n}&=\lim_{n\rightarrow \infty}
 \frac{b_{n+1}}{b_n}=17+12\sqrt{2}.
\end{align*}
With Mathematica,
 it is easy to verify
  that for $3\leq n\leq 1000$,
  \[
\frac{a_{n+1}}{a_n}<\frac{a_{n+2}}{a_{n+1}}.
  \]
For  the case $n\geq 1000$,
  it is easy to prove  that
 \begin{align}
h_1(n+1)-\frac{1}{(n+1)^2}> h_1(n)+
 \frac{1}{n^2}.\label{3-14}
 \end{align}
Thanks to
  \eqref{3-1} and \eqref{3-14},
 \begin{align}
\frac{a_{n+2}}{a_{n+1}}> h_1(n+1)-\frac{1}{(n+1)^2}
>h_1(n)+
 \frac{1}{n^2}>\frac{a_{n+1}}{a_n}.
 \end{align}
Therefore,
 the sequence  $\{a_{n+1}/a_n\}_{n\geq 3
 }$
 is    strictly increasing to the limit
   $17+12\sqrt{2}$.

   Using the same method, one can prove
    that the sequence  $\{
    b_{n+1}/b_n\}_{n\geq 1
 }$
 is    strictly increasing to the limit
   $17+12\sqrt{2}$. We omit the details.

  With Mathematica,
   it is easy to check
    that for $2 \leq n \leq 606$,
\[
\frac{\sqrt[n+1]{a_{n+1}}} {\sqrt[n]{a_n}}>\frac{\sqrt[n+2]{
a_{n+2}}}{\sqrt[n+1]{a_{n+1}}}
 .
\]
For the case $n\geq 607$,
  setting $(S_n,a_0,N_0,m,k_0)= (a_n,
17+12\sqrt{2}
 , 607,2,1)$
  and
  \[
v(n)=-\frac{\frac{153}{2}+54\sqrt{2}
 }{n}+\frac{\frac{675}{4}+\frac{7617\sqrt{2}}{64}
  }{n^2}
  \]
in Theorem \ref{Th-1}
 and using Lemma \ref{L-1},
  one can check that Conditions (1)-(4)
   are true. Therefore,  the sequence
 $\{\sqrt[n+1]{a_{n+1}}
 /\sqrt[n]{a_n}\}_{n\geq 607
  }$
  is  strictly    decreasing to 1.
   For the case $2 \leq n \leq 606$,

Furthermore, with Mathematica,
   one can verify
    that for $1 \leq n \leq 209$,
\[
\frac{\sqrt[n+1]{b_{n+1}}} {\sqrt[n]{b_n}}>\frac{\sqrt[n+2]{
b_{n+2}}}{\sqrt[n+1]{b_{n+1}}}
 .
\]
For the case $n\geq 210$,
  setting $(S_n,a_0,N_0,m,k_0)= (b_n,
17+12\sqrt{2}
 , 210,2,2)$
  and
  \[
v(n)=-\frac{\frac{85}{2}+30
 \sqrt{2}
 }{n}+\frac{\frac{275}{4}+\frac{3105
  \sqrt{2}}{64}
  }{n^2}
  \]
in Theorem \ref{Th-1}
 and using Lemma \ref{L-2},
  we  can check that Conditions (1)-(4)
   hold. Therefore,  the sequence
 $\{\sqrt[n+1]{b_{n+1}}
 /\sqrt[n]{b_n}\}_{n\geq 210
  }$
  is  strictly    decreasing to 1.
   This completes the proof.

\section{An analytic proof
of   Sun's conjecture}

   In this section, we present an
     analytic proof of Conjecture
      \ref{conjecture}
  based on  the following theorem
   proved by Xia \cite{Xia}.

\begin{theorem} \cite{Xia} \label{Th-2}
Let $\{S_n\}_{n\geq 1}$
 be a positive sequence. If there exist
  positive integers $k_0, N_0$
   and a function $f(n)$ such that
    $k_0<N_0^2+N_0+2$
     and for $n\geq N_0$,

     (i) $0<f(n)<\frac{S_n}{S_{n-1}}<f(n+1)
      $;

      (ii) $\frac{f(n+1)}{f(n+3)}>1-\frac{k_0}{
      n^2+n+2}$;

       (iii) $\left(1-\frac{k_0}{N_0^2
       +N_0+2}\right)^{N_0^2+N_0+2}f^{2N_0}(
        N_0)>S_{N_0}^2$;\\[6pt]
then for $n\geq N_0$,
\[
\frac{\sqrt[n+1]{S_{n+1}}}
{\sqrt[n]{S_n}}>\frac{\sqrt[n+2]{S_{n+2}}}{\sqrt[n+1]{S_{n+1}}}.
\]
\end{theorem}

In order
 to prove
 Conjecture
  \ref{conjecture},
   we require the following two lemmas.
    In fact, the   two
     functions $f_1(n)$
      and $f_2(n)$ given in the following
       two lemmas are constructed by the
       a heuristic approach
        stated in \cite{Xia}.

\begin{lemma}\label{L-3}
For $n\geq 5$,
\begin{align}\label{4-2-1}
f_1(n)<\frac{a_n}{a_{n-1}}<f_1(n+1),
\end{align}
where
\begin{align}\label{4-2-2}
f_1(n)=-\frac{(17+12\sqrt{2})
 (-256n^3+2304n^2-8352n+16245+444\sqrt{2}n
-3108\sqrt{2})}{32(2n-3)^3}.
\end{align}
\end{lemma}

\noindent{\it Proof.} We are ready to
  prove this lemma by induction on $n$.
 It is easy to check that
 \eqref{4-2-1} is true when
  $5\leq n\leq 8$.
   Assume that  \eqref{4-2-1} holds when $n=m-1$
    and $n=m$ with $m\geq 8$, namely,
 \begin{align}\label{4-2-4}
f_1(m-1)<\frac{a_{m-1}}{a_{m-2}}<f_1(m),
  \qquad
f_1(m)<\frac{a_{m}}{a_{m-1}}<f_1(m+1).
 \end{align}
We need  to prove
 that  \eqref{4-2-1} is true when $n=m+1
 $. Note that for
 $m\geq 8$,
  \begin{align}\label{4-2-5}
v_1(m)>0,\qquad v_2(m)<0,\qquad v_3(m)>0.
 \end{align}
It follows from \eqref{4-2-4} and \eqref{4-2-5}
  that
\begin{align}\label{4-2-6}
\frac{v_2(m+1)}{f_1(m)}<\frac{v_2(m+1)a_{m-1}
 }{a_m}<\frac{v_2(m+1)}{f_1(m+1)}.
\end{align}
  In view of \eqref{3-3} and \eqref{4-2-4}--\eqref{4-2-6},
\begin{align*}
\frac{a_{m+1}}{a_m}-f_1(m+1) = &v_1(m+1)
+v_2(m+1)\frac{a_{m-1}}{a_m} +v_3(m+1)\frac{a_{m-2}}{a_{m-1}}\cdot
 \frac{a_{m-1}}{a_m} -f_1(m+1)
 \nonumber\\[6pt]
 >&v_1(m+1)+\frac{v_2(m+1)}{
 f_1(m)}+\frac{v_3(m+1)
 }{f_1(m)f_1(m+1)}-f_1(m+1)
 \nonumber\\[6pt]
 =&\frac{3( 433-312\sqrt{2})
  t_1(m)}{230368 t_2(m)},
\end{align*}
where
\begin{align*}
t_1(m)=&-362337533952m^{11}
 +704128221184\sqrt{ 2}m^{10}
  +5395829489664m^{10}-36820724219904m^9
   \nonumber\\[6pt]
   &-12891237613568\sqrt{2}m^9
    +93584348512256\sqrt{2}m^8
     +167728449957888m^8
      \nonumber\\[6pt]
   &-372329602193024\sqrt{2}
       m^7-551533244000256m^7
    +897725199301760\sqrt{2}m^6
              \nonumber\\[6pt]
   &+1257420204922848
          m^6-1690830347739168m^5
    -1216174211891420\sqrt{2}m^5
                  \nonumber\\[6pt]
   &  +508645813354545m^4
             +369528698815820\sqrt{2}m^4
    +1925414020667475m^3
                     \nonumber\\[6pt]
   &  +1362363531867532\sqrt{2}m^3
                -1747022411934268\sqrt{2}
    -2466902398311021m^2
                                   \nonumber\\[6pt]
   & -1577508850402720\sqrt{2}m-
                   2231619479860815m
+48905330688-37342642176\sqrt{2}
\end{align*}
and
\begin{align*}
t_2(m)=&(-256m^3+2304m^2
 -8352m+444\sqrt{2}m+16245-3108\sqrt{2})
 \nonumber\\[6pt]
&\times (-256m^3+1536m^4-2-4512m+444\sqrt{2}m
+9941-2664\sqrt{2})(2m-1)^3(m+1)^3m.
\end{align*}
Note that for $m\geq 8$,
\[
\frac{3( 433-312\sqrt{2})
  t_1(m)}{230368 t_2(m)}>0
\]
and therefore,
\begin{align}\label{4-2-7}
\frac{a_{m+1}}{a_m}>f_1(m+1).
\end{align}
On the other hand, thanks
 to \eqref{3-3}
  and \eqref{4-2-4}--\eqref{4-2-6},
 \begin{align*}
\frac{a_{m+1}}{a_m}-f_1(m+2) = &v_1(m+1)
+v_2(m+1)\frac{a_{m-1}}{a_m} +v_3(m+1)\frac{a_{m-2}}{a_{m-1}}\cdot
 \frac{a_{m-1}}{a_m} -f_1(m+2)
 \nonumber\\[6pt]
 <&v_1(m+1)+\frac{v_2(m+1)}{
 f_1(m+1)}+\frac{v_3(m+1)
 }{f_1(m-1)f_1(m)}-f_1(m+2)
 \nonumber\\[6pt]
 =&\frac{3( 433-312\sqrt{2})
  r_1(m)}{230368 r_2(m)},
\end{align*}
where
\begin{align*}
r_1(m)=&-185516817383424m^{15}
 +4773675760877568m^{14}
  +447100189933568\sqrt{2}m^{14}
  \nonumber\\[6pt]
& -58453058434105344m^13 -8347711606620160\sqrt{2}m^{13}
 +80575696971235328\sqrt{2}m^{12}
   \nonumber\\[6pt]
&+419676942558560256m^{12}
 -1904280499445563392m^{11}
  -476538905346375680m^{11}
    \nonumber\\[6pt]
&+5481248449364508672m^{10}
 +1696505448326758400\sqrt{2}m^{10}
  -3090582218844547072\sqrt{2}m^9
    \nonumber\\[6pt]
&-8472149837116692480m^9-1529023473651028992
 \sqrt{2}m^8-843346900089828864m^8
   \nonumber\\[6pt]
&+35417136924482525568m^7+22673048553252742656
\sqrt{2}m^7-75062456368991906112m^6
  \nonumber\\[6pt]
&-50483727152054139424\sqrt{2}m^6
 +33943755572507411784\sqrt{2}m^5
  +47597507656839385290m^5
    \nonumber\\[6pt]
&+51453910934767700516\sqrt{2}m^4
 +74700627494547974223m^4
  -98795982089022044932\sqrt{2}m^3
    \nonumber\\[6pt]
&-139857386658156050895m^3 +2073072657457016637m^2
+1769483815471845708\sqrt{2}m^2
  \nonumber\\[6pt]
&+86104770356767537425m
 +60868524728579605260\sqrt{2}m
  +19452609284382720
     \nonumber\\[6pt]
& -13769080029642752\sqrt{2}
\end{align*}
and
\begin{align*}
r_2(m)=&(2m-3)(-256m^3+2304m^4-2-8352m
+444\sqrt{2}m+16245-3108\sqrt{2})\nonumber\\[6pt]
&\times (-256m^3+3072m^4-2-13728m+444\sqrt{2}m
 +27157-3552\sqrt{2})\nonumber\\[6pt]
 &\times (-256m^3+1536m^2
  -4512m+444\sqrt{2}m+9941-2664\sqrt{2})
   (2m+1)^3(m+1)^3m.
\end{align*}
It is easy to check that for $m\geq 8$,
\[
\frac{3( 433-312\sqrt{2})
  r_1(m)}{230368 r_2(m)}<0.
\]
Therefore,
\[
\frac{a_{m+1}}{a_m}-f_1(m+2)<0,
\]
from which with \eqref{4-2-7},
  we see that \eqref{4-2-1} holds
  when $n=m+1$. Thus, Lemma \ref{L-1}
    is proved
  by induction.   \qed

  By \eqref{3-11}, we can deduce the following
   lemma.

 \begin{lemma}\label{L-4}
   For $n\geq 4$,
   \begin{align}\label{4-3}
f_2(n)<\frac{b_n}{b_{n-1}}<f_2(n+1),
   \end{align}
where
\begin{align}\label{4-4}
f_2(n)=-\frac{(17+12\sqrt{2})
 (-768n^3+5376n^2-14304n
  +17599+180\sqrt{2}n-900\sqrt{2})}{96
   (2n-3)^3}.
\end{align}
  \end{lemma}

Since the proof of Lemma \ref{L-4} is analogous to that of Lemma
\ref{L-3},
  we omit the details.

To conclude this section,
 we turn to prove
  Conjecture \ref{conjecture}.

    By \eqref{4-2-1} and \eqref{4-3}, we see
that for $n\geq 8$,
 \begin{align}
\frac{a_{n+1}}{a_n}< &
 \frac{a_{n+2}}{a_{n+1}}, \label{4-5}
\\[6pt]
\frac{b_{n+1}}{b_n}< &
 \frac{b_{n+2}}{b_{n+1}}. \label{4-6}
 \end{align}
Moreover, one can check \eqref{4-5}
 is true
 when $3\leq n \leq 7$ and
  \eqref{4-6}
 holds
 when $1\leq n \leq 7$ with Mathematica.
  Therefore, in order to
   prove Conjecture \ref{conjecture},
    it suffices to prove that
    \begin{align}\label{4-7}
\frac{\sqrt[n+1]{a_{n+1}}} {\sqrt[n]{a_n}}>\frac{\sqrt[n+2]{
a_{n+2}}}{\sqrt[n+1]{a_{n+1}}} \qquad \qquad (n\geq 2)
\end{align}
and
  \begin{align}\label{4-8}
\frac{\sqrt[n+1]{b_{n+1}}} {\sqrt[n]{b_n}}>\frac{\sqrt[n+2]{
b_{n+2}}}{\sqrt[n+1]{b_{n+1}}}. \qquad \qquad (n\geq 1)
\end{align}

Taking $S_n=a_n$, $f(n)=f_1(n)$
 defined by \eqref{4-2-2}, $k_0=10$
  and $N_0=10$ in Theorem \ref{Th-2},
   we can check that the three conditions
    $(i)$, $(ii)$ and $(iii)$ hold.
     Hence, \eqref{4-7} is true for $n\geq 10$.
     With Maple, we can check
       \eqref{4-7} holds for  $2\leq n\leq 9$.

In addition, setting $S_n=b_n$, $f(n)=f_2(n)$
 defined by \eqref{4-4}, $k_0=5$
  and $N_0=5$ in Theorem \ref{Th-2},
   we  see  that the three conditions
    $(i)$, $(ii)$ and $(iii)$ are true
     and \eqref{4-8} holds
       for $n\geq 5$.
      For $1\leq n\leq 4$,
       we can verify  \eqref{4-8}  hold
        with Maple. This
        completes
         the proof. \qed

\noindent{\bf Acknowledgments.} We thank professor Zhi-Wei Sun for his conjectures and valuable suggestions. We also thank professor Qing-Hu Hou and Miss Han Wang for their helpful comments.  The first author is partially supported by the Natural Science Foundation of Jiangsu Province of China (\#BK20221383) and the National Natural Science Foundation of China (\# 11971203). The second author is supported by the National Natural Science Foundation of China(No.12001161)and Doctor Foundation of Heibei Normal University(No.L2019B05).

\end{document}